\documentclass[12pt, a4paper, reqno]{amsart}

\usepackage[utf8]{inputenc}
\usepackage{amsfonts,amssymb,amsxtra,amsthm,amsmath,amscd,mathrsfs,cite}
\usepackage[all]{xy}
\usepackage{enumerate}  
\usepackage{tikz}
\usepackage[normalem]{ulem}
\usepackage{soul}
\usepackage{color}

\newtheorem{theorem}{Theorem}[section]
\newtheorem{lemma}[theorem]{Lemma}
\newtheorem{proposition}[theorem]{Proposition}
	\theoremstyle{remark}
	
\numberwithin{equation}{section}
\addtocounter{footnote}{1}

\title[The M{\"o}bius Disjointness Conjecture on infinite-dimensional torus]
{The M{\"o}bius Disjointness Conjecture \\ on infinite-dimensional torus}

\author{Qingyang Liu, Jing Ma, Hongbo Wang}
\address{School of Mathematics, Renmin University of China, Beijing, 100872, People's Republic of China}
\email{qingyangliu@ruc.edu.cn}
\address{School of Mathematics, Jilin University, Changchun, 130012, People's Republic of China}
\email{jma@jlu.edu.cn} 
\address{School of Mathematics, Jilin University, Changchun, 130012, People's Republic of China}
\email{whb25@mails.jlu.edu.cn}
	

\begin{document}

\begin{abstract}
Let $\mathbb{T}^\omega$ be the infinite-dimensional torus, and  
$T: \mathbb{T}^\omega\to \mathbb{T}^\omega$ be defined by 
\[			
T: (x_1, x_2, \dots, x_k, \ldots)  
\mapsto 
(x_1 + \alpha, x_2 + h(x_1),  \dots, x_k + h(x_1 + (k-2)\beta), \dots)
\] 
with $\alpha\in \mathbb{R}, \beta\in \mathbb{R}\backslash\mathbb{Q},$ and 
$h: \mathbb{R}\to \mathbb{R}$ being $1$-period and $C^{1+\varepsilon}$-smooth.   
This flow $(\mathbb{T}^\omega, T)$ is distal, and is also irregular in the 
sense that its Birkhoff average does not exist for all $x\in \mathbb{T}^\omega$. 
The main result of this paper is that the M\"{o}bius Disjointness Conjecture of Sarnak holds for 
$(\mathbb{T}^\omega, T)$. 
\end{abstract}

\maketitle

\section{Introduction}
\subsection{The M\"{o}bius Disjointness Conjecture} 

Let $\mu(n)$ denote the M\"{o}bius function, defined by $\mu(1)=1$, $\mu(n)=(-1)^k$ if $n$ is a product of 
$k$ distinct primes, and $\mu(n)=0$ otherwise. 
The behavior of $\mu(n)$ is of fundamental importance in analytic number theory, 
particularly in the study of prime numbers. 
A flow is a pair $(X, T)$ where $X$ is a compact topological space and 
$T \colon X \to X$ a continuous map.  
Sarnak's M\"{o}bius Disjointness Conjecture \cite{Sar12} states the following.

\medskip

\noindent
\textbf{M\"{o}bius Disjointness Conjecture (Sarnak).}  
\emph{The M\"{o}bius function $\mu$ is linearly disjoint from every zero-entropy flow $(X,T)$.  
Equivalently, for any $f \in C(X)$ and any $x \in X$,
\begin{equation}\label{MDC}
\lim_{N\to\infty} \frac1N \sum_{n \leqslant N} \mu(n) \, f(T^n x) = 0 .
\end{equation}}

\medskip

This conjecture establishes a profound connection between number theory and 
dynamical systems, stimulating extensive research at the interface of analytic number theory and ergodic theory. 
The simplest instance occurs when $X$ is a single point and $T$ the identity map; this case is equivalent 
to the Prime Number Theorem.  
The first nontrivial case, where $X = \mathbb{T}$ is the unit circle and $T$ a translation, 
corresponds to Vinogradov's estimate for exponential sums over primes, which played a crucial role in 
his proof of the ternary Goldbach problem.  
To date, the conjecture has been verified for a growing class of flows; see for 
example \cite{Sar12, LiuSar15, FerKulLem18} and the references therein.

\subsection{Distal flows and skew products}

A flow $(X, T)$ is said to be distal if for any two distinct points $x \ne y$ in $X$,  
\[
\inf_{n \geqslant 0} d(T^n x, T^n y) > 0,
\]  
where $d$ is a metric on $X$. Results of Parry \cite{Par67} show that distal flows have zero topological entropy.  
Thus, in seeking further evidence for Sarnak's conjecture, a natural starting point is to examine the linear disjointness of $\mu$ from distal systems.  
According to Furstenberg's structure theorem for minimal distal flows \cite{Fur63}, skew products serve as the fundamental building blocks of distal flows. Consequently, the problem reduces to studying distal flows arising from skew-product constructions.

A skew product is a dynamical system defined on a product space with an action that couples the base and fiber components. 
A basic example is the skew product $T$ on the 2-torus $\mathbb{T}^2 = (\mathbb{R}/\mathbb{Z})^2$, given by  
\[
T: (x,y) \mapsto (x+\alpha, y + h(x)),
\]
where $\alpha \in \Bbb R$ and $h\in C(\Bbb R)$ is $1$-periodic.  
The flow $(\mathbb{T}^2, T)$ is distal, and its disjointness from the M\"obius function was first investigated by Liu and Sarnak \cite{LiuSar15, LiuSar17}. 
They proved that the M\"obius Disjointness Conjecture holds for $(\mathbb{T}^2, T)$ 
for all $\alpha \in \mathbb{R}$ and for analytic $h$ whose Fourier coefficients satisfy a certain 
technical condition. A notable feature of their result, especially from the perspective of KAM theory, is its validity for every $\alpha$, without Diophantine restrictions.  
Later, Wang \cite{Wan17} showed that analyticity of $h$ alone suffices, removing the extra Fourier condition.  
Subsequently, the regularity requirement on $h$ was further relaxed: Huang–Wang–Ye \cite{HuaWanYe19} proved 
the conjecture for $h$ of $C^\infty$-smooth, Kanigowski–Lemanczyk–Radziwill \cite{KanLamRad21} for 
$h$ of $C^{2+\varepsilon}$, and de Faveri \cite{Fav22} for $h$ of $C^{1+\varepsilon}$. 
We remark that if $h: \mathbb{R}\to \mathbb{R}$ is 
a $1$-period function of $C^r$-smooth for some $r>0$, then its Fourier coefficients satisfy 
$$
\hat{h}(n) \ll_r |n|^{-r}
$$ 
for all nonzero integers $n$. 

\subsection{Infinite-dimensional skew products} 
We define the infinite-dimensional torus ${\Bbb T}^\omega$
by the complete direct product ${\Bbb T}^\omega = \prod_{\Bbb N} {\Bbb T}$ 
of countably many copies of $\Bbb T$.  
The purpose of the paper is to establish the M\"obius Disjointness Conjecture for a wide class of nonlinear skew products  
$T: \mathbb{T}^{\omega} \to \mathbb{T}^{\omega}$ defined by  
\begin{align}\label{def/T}
T: 
\ &(x_1, x_2, \ldots x_k, \ldots) \nonumber\\
&\mapsto (x_1+\alpha, x_2+h(x_1), \ldots, x_k+h(x_1+(k-2)\beta), \ldots ),
\end{align}  
where $\alpha, \beta \in \mathbb{R}$ and $h$ is $1$-periodic continuous function.  
These flows $(\mathbb{T}^{\omega}, T)$ are distal by Auslander \cite{Aus88}, for which 
the M\"{o}bius Disjointness Conjecture is expected to hold. The main result of the paper is the following. 

\begin{theorem}[M\"obius Disjointness]\label{thm1}
Let $T$ be defined as in \eqref{def/T} with
$\alpha\in \mathbb{R}, \beta\in \mathbb{R}\backslash\mathbb{Q}$, 
and $h: \mathbb{R}\to \mathbb{R}$ being $1$-periodic function of $C^{1+\varepsilon}$-smooth.  
Then the M\"obius Disjointness Conjecture is true for $(\mathbb{T}^\omega, T)$.
\end{theorem}

Previously Liu \cite{QLiu19} extended the framework of Liu--Sarnak \cite{LiuSar17} to infinite-dimensional tori, 
establishing the M\"obius Disjointness Conjecture for certain skew products $(\mathbb{T}^{\omega}, Q)$  
named Furstenberg's irregular flows to be defined below. 
Specifically, let  
\begin{align}\label{def/Q} 
Q : (x_1, x_2, \ldots, x_k, \ldots) \mapsto 
( x_1+\alpha, x_2+h(x_1), \dots, \; x_k+h(x_1+\beta^{k-2}), \dots ),
\end{align}
where $\alpha$ is irrational with $\{q_m\}_{m\geqslant 1}$,  the sequence of denominators of its convergents,   
satisfying the exponential growth condition  
$q_{m+1} \asymp e^{q_m}$,  $\beta$ is irrational, and 
\begin{align}\label{Fur1}
h(x) = \sum_{m \neq 0} t_m (1 - e(q_m\alpha)) e(q_m x)
\end{align}
with all coefficients $t_m$ satisfying $|t_m| \leqslant \tau$ and $t_m = t_{-m}$ for some positive $\tau$. 
The flows $(\mathbb{T}^{\omega}, Q)$ are named Furstenberg's irregular flows, and 
Liu \cite{QLiu19} established the M\"obius Disjointness Conjecture for such flows. 

Note that Theorem \ref{thm1} holds for all $\alpha\in \mathbb{R}$, and also our 
function $h$ in Theorem \ref{thm1} is quite general, that is 
our $h$ is not required to satisfy the technical conditions related to \eqref{Fur1}. Thus 
Theorem \ref{thm1} contains M\"obius Disjointness for Furstenberg's irregular flows 
$(\mathbb{T}^{\omega}, Q)$ established in \cite{QLiu19}. 

The proof of Theorem \ref{thm1} naturally splits into two cases: whether $\alpha$ is rational or irrational. 
When $\alpha$ is rational, Theorem \ref{thm1} reduces to a classical theorem of Davenport \cite{Dav37} 
concerning exponential sums with the M\"{o}bius function. Thus, in this case the improvement originates 
from arithmetic properties.
When $\alpha$ is irrational, the improvement stems from dynamical systems theory. 
We establish in Theorem \ref{thm31} that the skew product $(\mathbb{T}^\omega, T)$ considered in Theorem \ref{thm1} has a polynomial rate of rigidity—a concept introduced by Kanigowski, Lemanczyk, and 
Radziwill \cite{KanLamRad21} that implies M{\"o}bius disjointness. Moreover, in 
Theorem \ref{thm32} we prove that under a slightly stronger hypothesis on $h$, the same 
skew product possesses sub-polynomial measure complexity—a notion due to Huang, Wang, 
and Ye \cite{HuaWanYe19} that likewise yields M{\"o}bius disjointness.
Since Theorems \ref{thm31} and \ref{thm32} are independent, we effectively provide two 
distinct proofs of Theorem \ref{thm1} under slightly different condition of $h$.

Finally we remark that the M\"obius Disjointness Conjecture on $\mathbb{T}^\omega$ 
in short intervals has recently been established in \cite{HeLiuMa}. 

\medskip 

\noindent\textbf{Notation.}  
We use standard notation in number theory. For example $e(x) = e^{2\pi i x}$, and $\|x\| = \min_{n \in \mathbb{Z}} |x-n|$.  
The notation $A \asymp B$ means $A \ll B$ and $B \ll A$.  
The set of positive integers is denoted by $\mathbb{N}$.

\section{Rigidity and  Measure complexity}\label{RM}
In this section, we collect some concepts and facts from \cite{KanLamRad21} and \cite{HuaWanYe19} without proof. 

Let $X$ be a compact space, $T:X\to X$ a homeomorphism and let $M(X, T)$ be the set of all $T$-invariant Borel probability measures on $X$.
	
\subsection{PR rigidity} For any $\nu\in M(X,T)$, the system $(X,T,\nu)$ is said to be rigid if for 
any $g\in L^2(X,\nu)$, there exists a strictly increasing sequence $\{r_n\}_{n\geqslant 1}\subseteq \mathbb{N}$ such that 
\[
\|g\circ T^{r_n}-g\|_{L^2(\nu)}^2 \to 0
\]
as $n\to \infty.$  
Moreover, $(X,T,\nu)$ is said to have polynomial rate of rigidity (PR rigidty) if there exists a linearly dense\footnote{A subset $\mathcal{F}$ of a space is linearly dense if the linear combinations of elements in $\mathcal{F}$ form  a dense subset.} subset $\mathcal{F}\subseteq C(X)$ such that for every $f\in \mathcal{F}$, there are $\delta>0$ and a strictly increasing sequence $\{r_n\}_{n\geqslant 1}\subseteq \mathbb{N}$ such that
\begin{eqnarray}\label{PRr}
\sum_{|j|\leqslant r_n^\delta} \|f\circ T^{jr_n}-f\|_{L^2(\nu)}^2 \to 0
\end{eqnarray}
as $n\to \infty.$   
The proposition below is a  criterion given in \cite[Theorem~1.1]{KanLamRad21}.
	
\begin{proposition}\label{prop2.1} 
If $(X,T, \nu)$ has PR rigidity for all $\nu\in M(X, T)$, then the M\"obius Disjointness Conjecture for $(X,T)$ holds. 
\end{proposition}
	
\subsection{Measure complexity} Suppose the topology of $X$ is induced by a metric $d$. For $n\in\mathbb{N}$, we define
\begin{equation}\label{bard}\notag
		\bar{d}_n(x, y):=\frac 1n\sum_{j=0}^{n-1}d(T^{j}(x), T^{j}(y))
\end{equation}
for $x, y\in X$. For $\nu \in M(X,T)$, the measure complexity of $(X,T)$ with respect to $\nu$ is concerned with the number of balls of radius $\varepsilon$ under the metric $\bar{d}_n$ that can cover the whole space $X$ except for a subset whose measure under $\nu$ is less than $\varepsilon$, for arbitrary $\varepsilon>0$. Below is the definition we need. 
	
Let $\nu \in M(X,T)$ and $\varepsilon>0$.  Define
\begin{align}\label{mc}
s_n(d, \nu, \varepsilon) 
& = s_n(X, T, d, \nu, \varepsilon)  \nonumber\\
&:= \min\bigg\{m\in\mathbb{N}: \nu\bigg(\bigcup_{j=1}^{m}B_{\bar{d}_n}(x_j, \varepsilon)\bigg) 
>1-\varepsilon \ \mbox{ for some $ x_1,\dots, x_m\in X$}\bigg\},
\end{align}
where $B_{\bar{d}_n}(x, \varepsilon)=\{y\in X:\bar{d}_n(x, y)<\varepsilon\}$. Note that $\nu(X)=1$. 
Let $\tau>0$ be fixed. We are going to require  
\begin{equation}\label{def/sub/pol}
\liminf_{n\rightarrow\infty}\frac{s_n(d, \nu, \varepsilon)}{n^\tau}=0
\end{equation}
for any $\varepsilon>0$. If \eqref{def/sub/pol} holds for all $\tau>0$, then we say 
that the measure complexity of $(X, T, \nu)$ is sub-polynomial. 

\medskip 
		
The following criterion given in \cite[Theorem 1.1, Propositions 2.2 and 4.1]{HuaWanYe19} 
provides a motivation for studying measure complexity.  

\begin{proposition}\label{mc}
If the measure complexity of $(X, T, \rho)$ is sub-polynomial for any $\rho \in M(X, T)$,
then the M{\"o}bius Disjointness Conjecture holds for $(X, T)$.
\end{proposition}
	
A flow $(Y,S,\nu)$ is said to be measurably isomorphic to $(X,T,\rho)$ if there are 
Borel measurable subsets $X'\subset X$ and $Y'\subset Y$ satisfying $\rho(X')=\nu(Y')=1$, $T X' \subset X'$, $S(Y') \subset Y'$ 
such that $\phi\circ T = S \circ \phi$ for some invertible measure-preserving map $\phi: X'\to Y'$.

\begin{proposition}\label{prop2.3}
Suppose $(Y,S,\nu)$ is measurably isomorphic to $(X,T,\rho)$. If the measure complexity of $(Y,S, \nu)$ is sub-polynomial, then so is $(X,T,\rho)$ and vice versa. 
\end{proposition}
	
\begin{proposition}\label{41}
Let $(X, T)$ be a flow and $\rho\in M(X, T)$. 
If $\rho$ has discrete spectrum, then the measure complexity of $(X, T, \rho)$ is a bounded function.
\end{proposition}
	
\section{Applications of dynamical criteria for disjointness}
 
The M\"obius Disjointness Conjecture has stimulated considerable interest in the dynamical systems community, 
leading to the development of several verifiable criteria.  

Let $T$ be a homeomorphism of a compact metric space $(X,d)$ and let $M(X,T)$ denote the set of $T$-invariant 
Borel probability measures on $X$. As we mentioned in the last section, 
Huang, Wang, and Ye \cite{HuaWanYe19} proved that if $(X, T, \nu, d)$ 
has sub-polynomial measure complexity 
for every $\nu\in M(X,T)$, then $(X,T)$ satisfies M\"obius disjointness. 
Kanigowski, Lemanczyk, and Radziwill \cite{KanLamRad21} showed that polynomial rate of rigidity for 
all $\nu\in M(X,T)$ also implies M\"obius disjointness.
Applying them, several authors have established M\"obius disjointness for various classes of systems, 
including flows on the $2$-torus $\mathbb{T}^2$ and on products $\mathbb{T} \times G/\Gamma$, 
where $G/\Gamma$ is a $3$-dimensional Heisenberg nilmanifold \cite{HuaLiuWan21,HeWan23, MaWu24, LauMa}.

As applications of these dynamical criteria, we are going to prove that the skew product defined 
in \eqref{def/T} is of sub-polynomial measure complexity, and also has polynomial rate of rigidity. 
In particular, we have the following theorems. 
 
\begin{theorem}[Polynomial rate rigidity]\label{thm31}
Let $T$ be defined as in \eqref{def/T} with 
$\alpha, \beta\in \mathbb{R}\backslash\mathbb{Q}$,  
and $h: \mathbb{R}\to \mathbb{R}$ being $1$-periodic of $C^{1+\varepsilon}$-smooth.  
Then the system $(\mathbb{T}^\omega, T, \nu)$ 
has a polynomial rate of rigidity for every invariant measure $\nu \in M(\mathbb{T}^\omega, T)$.
\end{theorem}

\begin{theorem}[Sub-polynomial measure complexity]\label{thm32}
Let $T$ be defined as in \eqref{def/T} with
$\alpha, \beta\in \mathbb{R}\backslash\mathbb{Q}$,  
and $h: \mathbb{R}\to \mathbb{R}$ being $1$-periodic of $C^{\infty}$-smooth. Then 
the measure complexity of $(\mathbb{T}^\omega, T, \nu)$ is 
sub-polynomial for every $\nu \in M(\mathbb{T}^\omega, T)$.
\end{theorem}

The M{\"o}bius disjointness of the flow $(\mathbb{T}^\omega, T)$ defined in \eqref{def/T} follows directly by applying 
Proposition \ref{prop2.1} to Theorem \ref{thm31} for $\alpha$ irrational. 
And the proof of Theorem \ref{thm1} for $\alpha$ rational will be given in Section \ref{alpharational}.

We note that Theorem \ref{thm1} can also be obtained, in the more restrictive 
case where $h$ is of $C^{\infty}$-smooth, by applying Proposition~\ref{mc} to Theorem~\ref{thm32}. 
This provides an independent proof of a restricted version of the main result for $\alpha$ irrational. 
Whether polynomial rate rigidity implies sub-polynomial measure complexity remains an open question.

The bulk of the paper is devoted to proving Theorems \ref{thm31} and \ref{thm32}.

\section{Preliminaries}
This section collects notation, facts, and auxiliary results that will be used throughout the paper.
\subsection{Continued fraction expansion of irrational $\alpha$}
Without loss of generality we assume that $\alpha \in [0,1)$. 
Such an irrational has its continued fraction expansion by
\[
\alpha = [0; a_1, a_2, \ldots, a_k, \ldots]=\frac{1}{a_1 + \frac{1}{a_2+\frac{1}{a_3+\ldots}}}.
\]
Let $l_k/q_k = [0; a_1, a_2, \ldots, a_k]$ be the $k$th convergent of $\alpha$. 
Then the sequence $\{q_k\}_{k\geqslant0}$ 
is strictly increasing and satisfies the following standard properties:

(P1) we have $q_0=1$, $q_1 = [1/\alpha]$ and, for all $k\geqslant 1$, 
$$
q_{k+1}=a_{k+1}q_k+q_{k-1}; 
$$ 

(P2) we have, for all $k\geqslant 1$, 
$$
\frac1{2q_{k+1}} < \|q_k \alpha\| < \frac1{q_{k+1}}.
$$ 
The following lemmas are from \cite[Lemmas 3.2 and 3.3]{Fav22}, and 
will be used in the proof of Theorem~\ref{thm31}.

\begin{lemma}\label{l1}
For any irrational $\alpha$ and any $k \geqslant 2$,
\[
\sum_{0<|q|<q_k}\frac{1}{\|q\alpha\|^2} \asymp q_k^2.
\]
\end{lemma}

\begin{lemma}\label{l2}
For any irrational $\alpha$, $k\geqslant 1$ and $1\leqslant c \leqslant q_k$,
\[
\sum_{q_k \leqslant |q| < q_{k+1}}\frac{1}{q^2}
\min\bigg\{\frac{1}{\|q\alpha\|^2},c^2 \bigg\} 
\ll \frac{c}{q_k}.
\]
\end{lemma}

\subsection{Two auxiliary sets $E$ and $M$}

For a given $\tau>0$, define
\begin{align}\label{E}
E_{\tau}:=\{k\geqslant 2:   
q_{k+1}>q_k^{\frac{1}{\tau}+3}\} 
\end{align}
and
\begin{equation}\label{M}
M:=\bigcup_{k\in E}\{\pm m_k q_k : m_k=1,2,\dots,a_k   \},
\end{equation}
where $a_k$ is as in (P1). In the proof of Theorem~\ref{thm32}, we distinguish between the cases where $M$ is finite or infinite. 
The following lemmas are \cite[Propositions 6.3, 6.4, and 6.2]{HuaWanYe19}, and will be applied accordingly.

\begin{lemma}\label{tpsi}
Let $h: \mathbb{R}\to \mathbb{R}$ be $1$-periodic of $C^\infty$-smooth, 
and $M$ finite. Define
\[
\widetilde{\psi}(t):=\sum_{m\neq 0}\hat{h}(m)\frac{1}{e(m\alpha)-1}e(mt),
\]
where $\hat{h}(m)$ is the $m${\color{red}{-}}th Fourier coefficient of $h$. Then $\widetilde{\psi}\in C(\mathbb{T})$ and
\[
\widetilde{\psi}(t+\alpha)-\widetilde{\psi}(t)=h(t)-\hat{h}(0).
\]
\end{lemma}

\begin{lemma}\label{l5}
Let $h: \mathbb{R}\to \mathbb{R}$ be $1$-periodic of $C^\infty$-smooth, and $M$ infinite. Then for any $\tau>0$ and $k\in E_\tau$,
\[
\max_{x \in \mathbb{T}}|H_{q_k}(x)-q_k\hat{h}(0)|\ll q_k^{-(\frac{1}{\tau}+2)}.
\]
\end{lemma}

\begin{lemma}\label{psi} 
Let $h: \mathbb{R}\to \mathbb{R}$ be $1$-periodic of $C^\infty$-smooth. Then the series 
\[
\sum_{m\notin M\cup \{0\}}\hat{h}(m)\frac{1}{e(m\alpha)-1}e(mt)
\]
converges uniformly to a continuous function $\psi$.
\end{lemma}

\subsection{A number-theoretic lemma for rational $\alpha$}
For the rational case in Theorem~\ref{thm1}, we require the following estimate of Davenport 
\cite{Dav37} on exponential sums agianst the M\"obius function. 

\begin{lemma}\label{Dav} 
Let $P(n)$ be a real polynomial of degree $d>0$ and $0\leqslant a<q$. Then for any $A>0$,
\[
\sum_{\substack{0<n\leqslant N \\ n \equiv a\bmod q}}\mu(n)e (P(n)) 
\ll \frac{N}{\log^{A} N},
\]
where the implied constant depends on $d, A, \varepsilon,$ and $q$, but not on the coefficients of $P(n)$.
\end{lemma}

\subsection{Analysis on $\mathbb{T}^\omega$}\label{Ana/TT}  
We review some basic facts about Fourier analysis on the infinite-dimensional torus; for proofs, see
e.g. Rudin \cite{Rud62}. 
As stated before, the infinite-dimensional torus ${\Bbb T}^\omega$
is defined by the complete direct sum 
${\Bbb T}^\omega = \prod_{\Bbb N} {\Bbb T}$ 
of countably many copies of $\Bbb T$.  We know that $\Bbb T$ is a compact abelian group, and so is ${\Bbb T}^\omega$
with the product topology. The dual group of ${\Bbb T}^\omega$ is the direct sum 
${\Bbb Z}^\infty = \oplus_{\Bbb N} {\Bbb Z}$
of countably many copies of $\Bbb Z$, where ${\Bbb Z}$ is the discrete abelian group of the integers.
Each $x\in {\Bbb T}^\omega$ may be thought of as a string $x=(x_1, \dots, x_k, \ldots)$,
and each $\gamma\in {\Bbb Z}^\infty$ can be written as $(n_1, \ldots, n_k, \ldots)$ where only
finitely many of the $n_k$'s are nonzero. The latter fact is important.

For any $f\in L^1({\Bbb T}^\omega)$, its Fourier transform $\hat{f}$ is
\begin{eqnarray}\label{def/Hat/f}
\hat{f}(n_1, \ldots, n_k, \ldots)
=\int_{{\Bbb T}^\omega}f(x) e(-\langle x, \gamma\rangle)dx
=\int_{{\Bbb T}^\omega}f(x) e\bigg(-\sum_{k}x_k n_k\bigg)dx,
\end{eqnarray}
where only finitely many of the integers $n_k$'s are different from $0$, and the $x_k$'s are real numbers modulo $1$.
Thus the inversion formula has the form
\begin{eqnarray}\label{def/Hat/f}
f(x_1, \ldots, x_k, \ldots)=\sum \hat{f} (n_1, \ldots, n_k, \ldots) e\bigg(\sum_{k}x_k n_k\bigg).
\end{eqnarray}
It is also known that the set of trigonometric polynomials on ${\Bbb T}^\omega$
is dense in $C({\Bbb T}^\omega)$.

We also make a remark on integration over product spaces. Let $I$ be an infinite index set 
and $J\subseteq I$ a finite subset. For each $i\in I$, let $(\Omega_i, \mathcal F_i, \nu_i)$ be a probability space. 
Let $\nu_{I}$ be a Borel probability measure on the product $\prod_{i\in I}\Omega_i$. 
If $A_{J} \subseteq \prod_{i \in J}\Omega_i$ is open and $\nu_{J}$ is the corresponding marginal measure, then
\[
\nu_{I} \bigg(A_{J}\times \prod_{i \in I-J}\Omega_i\bigg)=\nu_{J}(A_{J}).
\]
The integral of a function $f$ on $\prod_{i\in I}\Omega_i$ with respect to $\nu_I$ is denoted by 
\[
\int_{\prod_{i\in I}\Omega_i}{f(x)\mathrm{d}\nu_{I}(x)}. 
\]
These facts will be used in later sections. 

\section{Proof of Theorem \ref{thm31}} 

Define a metric $d$ on $\mathbb{T}^\omega$ by
\begin{equation}\label{d}
d(x,y)=\sum_{k=1}^{\infty}\frac{1}{2^k}\|x_k-y_k\|,
\end{equation}
where $x_k$ and $y_k$ are the $k$-th components of $x$ and $y$ respectively, 
and $\|x\|$ denotes the distance from $x$ to the nearest integer. 

For any $\varepsilon>0$, we will show that there exist an unbounded 
sequence $\{r_n\}_{n\geqslant 1}\subseteq \mathbb{N}$ such that \eqref{PRr} is true. 
Note that the triangle inequality and the $T$-invariance of $\nu \in M(\mathbb{T}^\omega, T)$ imply
\begin{equation*}
\| f \circ T^{jr_n} - f \|_{L^2(\nu)}^2 \leqslant \bigg(\sum_{i=1}^{|j|} 
\| f \circ T^{ir_n} - f \circ T^{(i-1)r_n} \|_{L^2(\nu)}\bigg)^2 
\leqslant j^2 \cdot \| f \circ T^{r_n} - f \|_{L^2(\nu)}^2, 
\end{equation*}
and therefore  
\begin{eqnarray*} 
\sum_{|j|\leqslant r_n^\delta} \|f\circ T^{jr_n}-f\|_{L^2(\nu)}^2 \leqslant r_n^{3\delta}\cdot \left\| f \circ T^{r_n} - f \right\|_{L^2(\nu)}^2.
	\end{eqnarray*}
	Since Lipschitz functions are dense in $C(\mathbb{T})$, 
	it suffices to consider the family of Lipschitz continuous functions $\mathcal{F}$ on $\mathbb{T}^\omega$. 
	For any  $f\in \mathcal{F}$, $n\in \mathbb{N}$ and $\nu \in M(\mathbb{T}^\omega, T)$, we have
	\begin{eqnarray*}
		\|f\circ T^n-f\|_{L^2(\nu)}^2 
		\ll_f
		\int_{\mathbb{T}^\omega}  
		d(T^n(x), x)^2 \textup{d} \nu(x).
	\end{eqnarray*}
	Taking $\delta=\varepsilon/400$ and $\lambda=\varepsilon/100$.
	To prove \eqref{PRr}, we will show that there exist an unbounded sequence $\{r_n\}_{n\geqslant 1}\subseteq \mathbb{N}$ such that 
	\begin{eqnarray}\label{A1}
		\int_{\mathbb{T}^\omega} d (  T^{r_n}(x), x )^2 \textup{d} \nu(x) \ll_\varepsilon 
		r_n^{-\lambda}
	\end{eqnarray}
	for any $\varepsilon>0$, 
	which implies
	\begin{eqnarray*} 
		\sum_{|j|\leqslant r_n^\delta} \|f\circ T^{jr_n}-f\|_{L^2(\nu)}^2 \leqslant r_n^{3\delta-\lambda}\rightarrow 0
	\end{eqnarray*}
	as $n\rightarrow \infty$.
	Then, the PR rigidity of $(\mathbb{T}^\omega, T)$ follows.   

After $n$ iterations, we see that for any $x=(x_1,x_2,\dots, x_k,\dots) \in  \mathbb{T}^\omega$,
	\begin{equation}\label{Tn}
		T^{n}(x_1,x_2,\dots, x_k,\dots ) = (y_1(n),y_2(n),\dots ,y_k(n),\dots ):=y(n),
	\end{equation} 
where $y_1(n)=x_1+n\alpha$ and 
\begin{equation*}\label{ykn}
y_k(n)=x_k+\sum_{j=0}^{n-1}h(x_1 +(k-2)\beta +j\alpha)
\end{equation*}
for $k\geqslant 2.$ 
It follows that 
$$
	d(T^{r_n}x,x)
	=\frac{1}{2}\|r_n \alpha\|
	+
	\sum_{k=2}^{\infty}
	\frac{1}{2^k}
	\|S_{h,r_n}(x_1+(k-2)\beta)\|,
$$
where
	$$
	S_{h,n}(t)=\sum_{j=0}^{n-1}h(t+j\alpha).
	$$
	Let $K_{r_n,\lambda}=K$ be the minimal integer $k$ such that  $2^{-2k} < r_n^{-\lambda}$, where $\lambda$ is a parameter to be chosen latter. 
	Then
	$$
	d(T^{r_n}x,x)^2 
	\ll \|r_n\alpha\|^2
	+\sum_{k=2}^{K}  \frac{1}{2^{2k}} \|S_{h,r_n}(x_1+(k-2)\beta)\|^2+r_n^{-\lambda}.
	$$
Taking integration, we have 
\begin{equation}\label{i}
		\int_{\mathbb{T}^{\omega}}{d(T^{r_n}x,x)^2 \mathrm{d}\nu(x)} 
		\ll \|r_n\alpha\|^2
		+\sum_{k=2}^{K} \frac{1}{2^{2k}} \int_{\mathbb{T}^{\omega}}
		\|S_{h,r_n}(x_1+(k-2)\beta)\|^2 \mathrm{d} \nu(x)+r_n^{-\lambda}.
	\end{equation}
Observing that the upper bound of this integrand is independent of the coordinates other than $x_1$,  
we can rewrite the integral on the right-hand side as
	$$
	\int_{\mathbb{T}^{\omega}}\|S_{h,r_n}(x_1+(k-2)\beta)\|^2 \mathrm{d}\nu(x)
	=
	\int_{\mathbb{T}}
	{\|S_{h,r_n}(x_1+(k-2)\beta)\|^2 \mathrm{d}_{(\pi_{*}\nu)}(x_1)},
	$$
	where $\pi_{*}\nu$ denotes the Borel probability measure on $\mathbb{T}$ induced by $\nu$. 
	Define
	$R_\alpha:\mathbb{T}\rightarrow \mathbb{T}$ by
	$x \mapsto x+\alpha$.
Since $\nu$ is $T$-invariant, $\pi_{*}\nu$ is $R_\alpha $-invariant. However, $\alpha$ is irrational, so $(\mathbb{T},R_\alpha)$ is uniquely ergodic. Noting that $\pi_{*}\nu$ is $R_\alpha $-invariant for any irrational $\alpha$, we see that $\pi_{*}\nu$ is the 
unique Borel probability measure on $\mathbb{T}$, i.e. the Lebesgue measure on $\mathbb{T}$. 
	
	Writing the Fourier expansion of $h$ as
	\begin{equation}\label{hF}\notag
		h(t)=\sum_{q\in \mathbb{Z}}    c_q e(qt),
	\end{equation}
	we get 
$$
S_{h,r_n}(x_1+(k-2)\beta)
= c_0r_n  +  \sum_{q \neq 0} c_q e(qx_1+q(k-2)\beta)\frac{1-e(qr_n \alpha)}{1-e(q\alpha)}.
$$
Therefore, noticing $S_{h,r_n}-c_0r_n \in L^2(\mathbb{T})$ and 
using Parseval's Theorem, we get
\begin{align*}
&\int_{\mathbb{T}}{\|S_{h,r_n}(x_1+(k-2)\beta)\|^2\mathrm{d}_{(\pi_{*}\nu)}(x_1)} \\ 
&\quad \ll
		\|c_0r_n\|^2
		+\int_{\mathbb{T}}  
		\bigg|\sum_{q \ne 0}c_qe(qx_1+q(k-2)\beta) \frac{1-e(qr_n \alpha)}{1-e(q\alpha)}
		\bigg|^2 \mathrm{d}x_1 
		\\
&\quad \ll
\|c_0r_n\|^2+\sum_{q \ne 0}|c_q|^2 
\bigg|\frac{1-e(qr_n \alpha)}{1-e(q\alpha)}\bigg|^2. 
\end{align*}
The last expression is independent of $k$. So we rewrite \eqref{i} as 
	\begin{align}\label{i2}
		\int_{\mathbb{T}^{\omega}}{d(T^{r_n}x,x)^2 d\nu(x)} 
		\ll \|r_n\alpha\|^2
		+  \|c_0r_n\|^2
		+   
		\sum_{q \ne 0}|c_q|^2 
		\bigg|\frac{1-e(qr_n \alpha)}{1-e(q\alpha)}\bigg|^2+r_n^{-\lambda}.
	\end{align}

	Now we will select an appropriate strictly increasing sequence $\{r_n\}_{n=1}^{\infty}$.
	Let  $q_n$ be the denominator of the $n$-th convergent of the continued fraction expansion of $\alpha$.
	The Dirichlet approximation theorem implies that there is an $l_n \in \mathbb{Z}$ such that 
	$$
	1 \leqslant l_n \leqslant q_n^{\gamma}, \quad \|c_0l_nq_n\|<q_n^{-\gamma},
	$$
	where $\gamma :={\varepsilon}/{10}$. 
	Let $r_n:=l_nq_n$. Then $r_n\leqslant q_n^{1+\gamma}$. Recalling $\lambda= {\varepsilon}/{100}$, we get 
	\begin{align}\label{i21}
		\|r_n\alpha\|^2
		\leqslant  l_n^2\|q_n\alpha\|^2<q_n^{2\gamma}q_{n+1}^{-2}<q_n^{2\gamma-2}<q_n^{-\lambda(1+\gamma)}\leqslant  r_n^{- \lambda} 
	\end{align}
	by applying (P2), and 
	\begin{align}\label{i22}
		\|c_0r_n\|^2=\|c_0l_nq_n\|^2 <q_n^{-2\gamma}<q_n^{-\lambda(1+\gamma)} \leqslant r_n^{-\lambda}. 
	\end{align}
To deal with the sum over $q$ in \eqref{i2}, we split it into two parts according as $|q|\geqslant q_n$ or $0<|q|<q_n$. 
For the first part, we have 
\begin{align}
		\sum_{|q|\geqslant q_n}|c_q|^2 \bigg|\frac{1-e(qr_n\alpha)}{1-e(q\alpha)}\bigg|^2
		\ll 
		& \sum_{|q|\geqslant q_n}\frac{1}{|q|^{2+2\varepsilon}}\min\bigg\{\frac{1}{\|q\alpha\|^2},r_n^2 \bigg\}
		\notag 
		\\
		\ll &q_n^{-2\varepsilon}l_n^2\sum_{k=n}^{\infty}\sum_{q_k \leqslant |q| <q_{k+1}}\frac{1}{q^2} 
		\min \bigg\{ \frac{1}{\|q\alpha\|^2},q_n^2 \bigg\}\notag 
		\\
		\ll &q_n^{-2\varepsilon+2\gamma}\sum_{k=n}^{\infty}\frac{q_n}{q_k}
		\ll q_n^{-2\varepsilon+2\gamma}<q_n^{-\lambda(1+\gamma)}\leqslant r_n^{-\lambda},\nonumber
\end{align}
where we have used Lemma \ref{l2} (for $c=q_n \leqslant q_k$) and the fact that $q_{k+2} >2q_k$ by (P1) which implies that $q_k$ grow exponentially, together with 
$$
c_q\ll \frac{1}{|q|^{1+\varepsilon}}, \quad |1-e(q\alpha)| \asymp \|q\alpha\|. 
$$ 
For the second part with $0 <|q|<q_n$, we have 
$$
	|1-e(qr_n\alpha)|^2 \asymp \|ql_nq_n \alpha\|^2 \leqslant q^2l_n^2\|q_n\alpha\|^2 < q^2q_n^{2\gamma}q_{n+1}^{-2}, 
$$
and therefore 
\begin{equation}\label{scq}
		\sum_{0<|q|<q_n}|c_q|^2 \bigg|\frac{1-e(qr_n\alpha)}{1-e(q\alpha)}\bigg|^2
		\ll_h 
		\frac{q_n^{2\gamma}}{q_{n+1}^2}\sum_{0<|q|<q_n}\frac{1}{|q|^{2\varepsilon}} \frac{1}{\|q\alpha\|^2}.
\end{equation} 

To deal with this sum, we consider two cases. 

\medskip 
\noindent 
{\sc Case 1:} Suppose that $q_{n+1}<q_{n}^2$ for all sufficiently large $n$. Without loss of generality, we may assume $q_{n+1}<q_{n}^2$ for all $n\geqslant1 $. Take $0<k<n$ such that $q_k \in [q_n^{1/4},q_n^{1/2}]$, then $q_k^2\leqslant q_n$ and $q_k^{-1}\leqslant q_n^{-1/4}$. 
	Using Lemma \ref{l1}, we can rewrite the sum on the right-hand side of \eqref{scq} as
\[
\bigg\{\sum_{0<|q|<q_k}+\sum_{q_k \leqslant |q| <q_n}\bigg\}
\frac{1}{|q|^{2\varepsilon}} \frac{1}{\|q\alpha\|^2}<\sum_{0<|q|<q_k}\frac{1}{\|q\alpha\|^2}+q_k^{-2\varepsilon}\sum_{q_k \leqslant |q| <q_n}\frac{1}{\|q\alpha\|^2} \ll q_k^2+q_k^{-2\varepsilon}q_n^2.
\]	
Then the corresponding upper bound of \eqref{scq} is
\begin{equation}\label{i23}
		\sum_{0<|q|<q_n}|c_q|^2 \bigg|\frac{1-e(qr_n\alpha)}{1-e(q\alpha)}\bigg|^2
		\ll
		\frac{q_n^{2\gamma}}{q_{n+1}^2} (q_n+q_n^{2-\varepsilon/2})
		\ll q_n^{2\gamma-\varepsilon/2}
		<q_n^{-\lambda(1+\gamma)}\leqslant r_n^{-\lambda}.
\end{equation}  
	
\medskip 
\noindent 
{\sc Case 2:} Suppose that there is a subsequence $\{q_{b_n}\}_{n \geqslant 1}$ of $\{q_n\}_{n \geqslant 1}$ such that $q_{b_{n+1}}\geqslant q_{b_n}^2$ for all $n \geqslant 1$.
	Using \eqref{scq}, $q_{b_{n+1}}^{-2}\leqslant q_{b_n}^{-4}$ and Lemma \ref{l1} we get an upper bound for \eqref{scq} as
	\begin{align}
		\sum_{0<|q|<q_{b_n}}|c_q|^2 \bigg|\frac{1-e(qr_{b_n}\alpha)}{1-e(q\alpha)}\bigg|^2
		\ll  
		\frac{q_{b_n}^{2\gamma}}{q_{b_{n}}^2}
		\sum_{0<|q|<q_{b_n}} \frac{1}{q_{b_n}^2}\frac{1}{\|q\alpha\|^2} 
		\ll 
		q_{b_{n}}^{2\gamma-2}<q_{b_n}^{- \lambda(1+\gamma)} \leqslant r_{b_n}^{-\lambda}.\nonumber
	\end{align}
	By taking the subsequence $\{r_{b_n}\}_{n\geqslant 1}$ instead of the original sequence  $\{r_{n}\}_{n\geqslant 1}$ we get \eqref{i23} again.   
	
Summarizing \eqref{i2}, \eqref{i21}, \eqref{i22} and \eqref{i23}, 
we get \eqref{A1} and hence Theorem~\ref{thm31} for irrational $\alpha$. 
	
\section{Proof of Theorem \ref{thm1}}\label{alpharational}

As pointed out before, Theorem \ref{thm1} for irrational $\alpha$ can be deduced from Theorem \ref{thm31}, 
hence we are left with the case that $\alpha$ is rational. 
Let $\alpha={l}/{q}$ for some $l,q\in \mathbb{Z}$. 
	After $n$ iterations, we see that for any $x \in \mathbb{T}^\omega$,
$T^{n}(x)$ is of the form  \eqref{Tn}.
	Since the space of trigonometric functions is dense in $C(\mathbb{T}^\omega)$, 
	it is sufficient to show that for a given  
	$b\in\mathbb{Z}^\infty $
	we have 
	$$
	S(N)
	:=\sum_{1\leqslant n\leqslant N}\mu(n) e (\langle b, T^n(x)\rangle)
	\ll N(\log N)^{-A},
	\qquad
	\forall
	x\in\mathbb{T}^\omega,
	$$
	for some $A>0$, 
	where $\langle b, T^n(x)\rangle$ means the dot product.   

By the structure of $\mathbb{Z}^\infty$ in \S\ref{Ana/TT}, any $b\in \mathbb{Z}^\infty$ must be of the form 
$b=(b_1,\dots,  b_{K},0,\dots )$. It follows from this and \eqref{Tn} that 
\[
\langle b,T^n(x)\rangle
	=b_1(x_1+n\alpha)
	+\sum_{k=2}^{K}b_k x_k
	+\sum_{k=2}^{K}b_{k}\sum_{j=0}^{n-1}h (x_1 +(k-2)\beta +j\alpha). 
\]
	Then, we can write
	\begin{align}
		\label{smn}
		S(N)
		&=e(\sum_{k=1}^{K}b_k x_k )\sum_{1\leqslant n\leqslant N}\mu(n)e\bigg\{b_1n\alpha+\sum_{k=2}^{K}b_k \sum_{j=0}^{n-1}h (x_1 +(k-2)\beta +j\alpha)\bigg\}\notag
		\\
		&\ll\sum_{1\leqslant n\leqslant N}\mu(n)e\bigg\{b_1n\alpha+\sum_{k=2}^{K}b_k\sum_{j=0}^{n-1}h (x_1 +(k-2)\beta +j\alpha)\bigg\}.
	\end{align}

	If $\alpha=0$, using Lemma \ref{Dav}, we can write  \eqref{smn} as
	\[
	S(N)
	\ll
	\sum_{1\leqslant n\leqslant N}\mu(n)e\bigg\{n\sum_{k=2}^{K}b_k h (x_1+(k-2)\beta)\bigg\}
	\ll_{A,\varepsilon}
	\frac{N}{\log^A N}.
	\]
	
	Suppose that $\alpha\neq 0$.  Divide the last sum over $j$ in \eqref{smn} into residue classes modulo $q$. 
	Since $h(x)=h(x+1)$ for any $x\in\mathbb{T}$, 
	letting $n\equiv r\bmod{q}$ for some $0\leqslant r\leqslant q-1$, and
	$$
	\gamma_1^{(r)}:=\sum_{j=0}^{r-1}h\bigg(x_1+(k-2)\beta+j\frac{l}{q}\bigg),\quad
	\gamma_2^{(r)}:=\sum_{j=r}^{q-1}h\bigg(x_1+(k-2)\beta+j\frac{l}{q}\bigg),
	$$
	using Lemma \ref{Dav}, we get
\begin{align*}
S(N) 
&\ll \sum_{r=0}^{q-1}\sum_{\substack{n\equiv r \bmod q\\ 1\leqslant n\leqslant N}}
\mu(n)e\bigg( b_1n\alpha +\sum_{k=2}^{K} b_k \bigg(\gamma_1^{(r)}\frac{n-r+q}{q}+\gamma_2^{(r)}\frac{n-r}{q}\bigg)\bigg)\\ 
&\ll_{A,q,\varepsilon} \frac{N}{\log^{A}N}.
\end{align*}
This completes the rational case of Theorem~\ref{thm1}.  

\section{Proof of Theorem \ref{thm32}}
Since $\alpha$ is irrational, we can define $E$ and $M$ as 
in \eqref{E} and \eqref{M} respectively. We will prove Theorem \ref{thm32} in two cases according to 
the set $M$ being finite or infinite.
	
\subsection{$M$ is finite}
Define two maps $\widetilde{S}$ and $\widetilde{\pi}: \mathbb{T}^\omega\rightarrow \mathbb{T}^\omega$ by
\begin{align*}
\widetilde{S}: (x_1, x_2, \ldots, x_k,\ldots) 
\mapsto
(x_1+\alpha, x_2+\hat{h}(0), \ldots, x_k+\hat{h}(0),\ldots),
\end{align*}
and 
\begin{align*}
\widetilde{\pi}: (x_1, x_2, \ldots, x_k,\ldots) 
\mapsto
(x_1,x_2-\widetilde{\psi}(x_1),\dots,x_k-\widetilde{\psi}(x_1+(k-2)\beta),\dots),
\end{align*}
where $\widetilde{\psi}$ is defined in Lemma \ref{tpsi}.
Since $h: \mathbb{R}\to \mathbb{R}$ is $1$-periodic of $C^\infty$-smooth, 
and $M$ is finite, we can  apply Lemma \ref{tpsi} to get
\begin{align*}
&\widetilde{\pi}^{-1}  \circ \widetilde{S} \circ \widetilde{\pi}  (x_1, x_2, \ldots, x_k,\ldots)  
\\
&=
\widetilde{\pi}^{-1}  \circ \widetilde{S} (x_1,x_2-\widetilde{\psi}(x_1),\dots,x_k-\widetilde{\psi}(x_1+(k-2)\beta),\dots )
\\
&=
\widetilde{\pi}^{-1}  
(x_1+\alpha, x_2-\widetilde{\psi}(x_1) +\hat{h}(0),
\dots,
x_k-\widetilde{\psi}(x_1+(k-2)\beta) +\hat{h}(0),\dots)
\\
&=
(x_1+\alpha, x_2+\widetilde{\psi}(x_1+\alpha)  -\widetilde{\psi}(x_1) +\hat{h}(0),
\dots,
\\
& \qquad \qquad 
		x_k+\widetilde{\psi}(x_1+\alpha+  (k-2)\beta)  -\widetilde{\psi}(x_1+(k-2)\beta) +\hat{h}(0),\dots). 
\end{align*} 
The last line can be simplified as 
\begin{align*}
&= (x_1+\alpha, x_2+h(x_1), \dots, x_k+h(x_1 +(k-2)\beta),\dots) \\ 
&= T(x_1, x_2, \dots,x_k,\dots),  
\end{align*} 
that is $T$ is conjugate to $\widetilde{S}$ by $\widetilde{\pi} $. 

Let $\widetilde{\nu}= \rho \circ \widetilde{\pi}^{-1}$. 
Then $(\mathbb{T}^{\omega},T, \rho)$ is measurably isomorphic to $(\mathbb{T}^\omega,\widetilde{S},\widetilde{\nu})$. 
Since $\widetilde{S}$ is a rotation of $\mathbb{T}^\omega$, 
all the invariant measures of $T$ have discrete spectrum and bounded measure complexity by Proposition \ref{41}, 
and therefore the measure complexity of $(\mathbb{T}^\omega,\widetilde{S},\widetilde{\nu})$ is sub-polynomial.
Hence the measure complexity of $(\mathbb{T}^{\omega},T, \rho)$ is also sub-polynomial.
	
\subsection{$M$ is infinite}
Now we consider the case where $M$ is infinite. 
Define
\begin{equation}\label{h1}
h_1(t):=\sum_{m\in M \cup \{0\}}\hat{h}(m)e(mt).
\end{equation}
Then $h_1\in C^\infty$ and 
$$
h(t)-h_1(t)= {\psi}(t+\alpha)  - {\psi}(t)
$$ 
for any $t\in [0, 1)$, where $\psi$ is defined in Lemma \ref{psi}.
Define two maps $S$ and $\pi: \mathbb{T}^\omega\rightarrow \mathbb{T}^\omega$ by
\begin{align*}
&    S: (x_1, x_2, \dots, x_k, \ldots)
\mapsto 
(x_1+\alpha, x_2+h_1(x_1), 
		\dots,x_k+h_1(x_1+(k-2)\beta),\dots),
\end{align*}  
and 
\begin{align*}
\pi: (x_1, x_2, \dots, x_k, \ldots)  
\mapsto 
(x_1,x_2-\psi(x_1),  \dots,x_k-\psi(x_1+(k-2)\beta),\dots).
\end{align*}  
Then,  
\begin{align*}
		&{\pi}^{-1}  \circ S  \circ {\pi}  (x_1, \dots,x_k,\dots) 
		\\
		&=
		{\pi}^{-1}  \circ S  (x_1,x_2-{\psi}(x_1),\dots,x_k-{\psi}(x_1+(k-2)\beta),\dots)
		\\
		&=
		{\pi}^{-1}  
		(x_1+\alpha, x_2-{\psi}(x_1) + {h}_1(x_1),
		\dots,
		x_k-{\psi}(x_1+(k-2)\beta) + {h}_1(x_1+(k-2)\beta),\dots)
		\\
		&=
		(x_1+\alpha, x_2+{\psi}(x_1+\alpha)  - {\psi}(x_1) + {h}_1(x_1),
		\dots,
		\\
& \qquad\qquad x_k+ {\psi}(x_1+\alpha  +(k-2)\beta)  - {\psi}(x_1+(k-2)\beta) + {h}_1(x_1+(k-2)\beta),\dots). 
\end{align*} 
The last line is 
\begin{align*}
&=
(x_1+\alpha, x_2+h(x_1), \dots, x_k+h(x_1 +(k-2)\beta),\dots)\\
&=T(x_1, \dots,x_k,\dots). 
\end{align*} 
This means that $T$ is conjugate to $S$ by $\pi$. 
Let $\nu = \rho \circ \pi^{-1}$. 
Then $(\mathbb{T}^\omega,{S},{\nu})$ is measurably isomorphic to $(\mathbb{T}^{\omega},T, \rho)$. 
By Proposition \ref{prop2.3}, to prove that the measure complexity of $(\mathbb{T}^{\omega},T, \rho)$ is sub-polynomial, 
it will be enough to prove that the measure complexity of $(\mathbb{T}^\omega,S,\nu)$ is sub-polynomial. 

\subsection{Measure complexity of $(\mathbb{T}^\omega,S,\nu)$}  
In this subsection we are going to prove that the measure complexity of $(\mathbb{T}^\omega,S,\nu)$ is sub-polynomial.  
Let $\tau>0$ be arbitrarily small. Let $s_n(d,\nu,\varepsilon)$ be defined in \eqref{mc}, 
and the metric $d$ be as in \eqref{d}. We want to prove that 
\begin{equation}\label{Claim/lim}
\lim_{n\rightarrow +\infty} \inf\frac{s_n(d,\nu,\varepsilon)}{n^{\tau}}=0
\end{equation} 
holds for any $\varepsilon>0$. The rest of this section is devoted to the proof of  \eqref{Claim/lim}. 

We define
$$
H_n(t):=\sum_{j=0}^{n-1}h_1(t+j\alpha)
$$
for $n\in \mathbb{N}$ and $H_0(x)=0$. 
Then, after $n$ iterations, we get
$$
S^n(x_1, x_2, \dots,x_k,\dots) 
=
(x_1+n\alpha, x_k+H_n(x_1), \dots, x_k+H_n(x_1+(k-2)\beta),\dots).
$$

Since $M$ is infinite, for any $\varepsilon >0$, by Lemma \ref{l5}, there exists a positive constant $C$ such that
$$
\max_{x \in \mathbb{T}}|H_{q_t}(x)-q_t\hat{h}(0)|\leqslant  C q_t^{-(\frac{1}{\tau}+2)}
$$
for all $t\in E$, where $E$ is defined in \eqref{E}. 
This implies that
\begin{align}\label{Hqt}
\max_{x_i, x_i^* \in \mathbb{T}}|H_{q_t}(x_i)-H_{q_t}(x_i^*)| \leqslant  2C q_t^{-(\frac{1}{\tau}+2)}
\end{align}
for all $t \in E$.
Choose $t_0\in \mathbb{N}$ such that 
\begin{equation}\label{t}
\frac{2C}{q_t}<\frac{\varepsilon}{4}
\end{equation} 
for all $t \geqslant t_0$. For $t\in E\cap [t_0,+\infty)$, let 
\begin{align}\label{Hq_t}
		n_t :=q_t^{[\frac{1}{\tau}]+2},
	\end{align}
	where $[\frac{1}{\tau}]$ is the integer part of $\frac{1}{\tau}$.
	Since $h_1(x)\in C^\infty$, 
	there exists an integer $L\geqslant {4}/{\varepsilon}$ such that
	\begin{align}\label{Lx_2}
		|h_1(x_i)-h_1(x_i^*)| \leqslant  L \|x_i-x_i^* \| 
	\end{align}
	for any $x_i,x_i^*\in \mathbb{T}$.
	
	For each $\varepsilon>0$, there exists an $N\in\mathbb{N}$ such that
	$$
	\sum_{n\geqslant N+1}\frac{1}{2^n}<
	\frac{\varepsilon}{4}.
	$$
	For $t\in E\cap [t_0,+\infty)$,
let
\begin{align*}
F_t: 
=\bigg\{\bigg(
\frac{i}{Lq_t([\frac{4}{\varepsilon}]+1)}, \frac{j_2}{L},\dots,\frac{j_N}{L},0,\dots \bigg): \ 
& 0  \leqslant  i  \leqslant Lq_t\bigg(\bigg[\frac{4}{\varepsilon}\bigg]+1\bigg) -1, \\
& 0  \leqslant  j_k \leqslant  L-1, \ 2 \leqslant k \leqslant N  \bigg\}.
\end{align*} 
	For any $x=(x_1,x_2,\dots)\in \mathbb{T}^\omega$, there is an $x^*=(x_1^*,x_2^*,\dots) \in F_t$ such that
\begin{equation}\label{ll}
\|x_1-x_1^* \| 
\leqslant  
\frac{1}{Lq_t([\frac{4}{\varepsilon}]+1)}
\leqslant\frac{1}{L}\leqslant \frac{\varepsilon}{4},
\quad \|x_k-x_k^* \|\leqslant  \frac{1}{L} \leqslant  \frac{\varepsilon}{4},
\end{equation} 
for $k=2, 3, \dots, N.$ For $i = 0, 1, \dots, n_t-1$, write
$$
i=a_iq_t+b_i
$$
with 
\begin{equation}\label{ab}
0   a_i \leqslant q_t^{[\frac{1}{\tau}]+1}-1,  
\quad
0\leqslant  b_i \leqslant  q_t-1.
\end{equation}  
Then we have 
\[
S^i(x)=\bigg(x_1+i\alpha, \dots,x_k+\sum_{r=0}^{a_i-1}H_{q_t}(x_1+(k-2)\beta+rq_t\alpha)+\sum_{j=0}^{b_i-1}h_1(x_1+(k-2)\beta+(a_iq_t+j)\alpha),\dots \bigg)
\]	
and 
\[	
S^i(x^*)=\bigg(x_1^*+i\alpha, \dots,x_k^*+\sum_{r=0}^{a_i-1}H_{q_t}(x_1^*+(k-2)\beta+rq_t\alpha)+\sum_{j=0}^{b_i-1}h_1(x_1^*+(k-2)\beta+(a_iq_t+j)\alpha),\dots \bigg). 
\]	
These together with \eqref{ll}, \eqref{ab}, \eqref{Hqt}, \eqref{Lx_2} and \eqref{t} imply that 
\begin{align}
d(S^i(x),S^i(x^*))
&\leqslant 
\frac{1}{L}+q_t^{[\frac{1}{\tau}]+1} \cdot 2Cq_t^{-(\frac{1}{\tau}+2)}+q_t \cdot  L \cdot \frac{1}{Lq_t([\frac{4}{\varepsilon}]+1)}
+\frac{\varepsilon}{4}
\notag\\
&\leqslant 
\frac{\varepsilon}{4}+\frac{2C}{q_t}+\frac{1}{[\frac{4}{\varepsilon}]+1}+\frac{\varepsilon}{4}
\leqslant 
\varepsilon.   
\end{align}
It follows that 
\[
\overline{d}_{n_t} (x,x^*)
=
\frac{1}{n_t}\sum_{i=0}^{n_t-1}d(S^i(x),S^i(x^*)) \leqslant \varepsilon,
\]
which means that 
\[ 
x \in B_{\overline{d}_{n_t}}(x^{*},\varepsilon) \in \bigcup_{x^*\in F_t}B_{\overline{d}_{n_t}}(x^*,\varepsilon). 
\] 
Therefore 
$\cup_{x^*\in F_t}B_{\overline{d}_{n_t}}(x^*,\varepsilon)=\mathbb{T}^\omega.$ 
Hence 
	\begin{align}\label{SFt}
		s_{n_t}(d,\nu,\varepsilon) \leqslant |F_t|
		=L^{N} q_t \bigg(\bigg[\frac{4}{\varepsilon}\bigg]+1\bigg).
	\end{align}
	
Finally, since $M$ and $E$ are infinite, we have by \eqref{SFt} that 
\begin{align*}
		\lim_{n \rightarrow +\infty}\inf \frac{s_{n}(d,\nu,\varepsilon)}{n^\tau}
		&\leqslant  
		\lim_{\substack{t \rightarrow +\infty  \\
				t\in E\cap[t_0,+\infty)}}\inf \frac{s_{n_t}(d,\nu,\varepsilon)}{n_t^\tau} 
		\leqslant  \lim_{\substack{t \rightarrow +\infty \\ \notag
				t\in E\cap[t_0,+\infty)}}\inf \frac{L^{N}q_t([\frac{4}{\varepsilon}]+1)}{q_t^{\tau([\frac{1}{\tau}]+2)}}\\ \notag
&\leqslant
\lim_{\substack{t \rightarrow +\infty \\ \notag
t\in E\cap[t_0,+\infty)}}\inf \frac{L^{N}([\frac{4}{\varepsilon}]+1)}{q_t^{2\tau}}=0.
\end{align*}
This completes the proof of \eqref{Claim/lim} for any $\tau>0$, and hence the proof of Theorem \ref{thm32}. 
	
\subsection*{Funding} Liu is supported by the National Natural Science 
Foundation of China (Grant No. 12401012).


\end{document}